\documentclass[10pt]{article}
\usepackage{amssymb,amsmath,amsthm}
\newcommand{\ie}{\emph{i.e.}}
\newcommand{\eg}{\emph{e.g.}}
\newcommand{\cf}{\emph{cf}}

\newcommand{\demi}{\frac{1}{2}}
\newcommand{\Real}{\mathbb{R}}

\newcommand{\Nat}{\mathbb{N}}

\newcommand{\vol}{\mathrm{vol}}
\newcommand{\Smooth}{C}

\newcommand{\sii}{L^2}
\newcommand{\sinf}{L^\infty}
\newcommand{\Dom}{\mathop{\mathrm{Dom}}\nolimits}
\newcommand{\supp}{\mathop{\mathrm{supp}}\nolimits}

\newcommand{\diag}{\mathop{\mathrm{diag}}\nolimits}
\newcommand{\tr}{\mathop{\mathrm{tr}}\nolimits}

\newcommand{\sobi}{\mathop{W_0^{1,2}}\nolimits}
\newcommand{\Sobi}{\mathop{W^{1,2}}\nolimits}

\newcommand{\Id}{1}
\newcommand{\ef}{\mathcal{J}}
\newcommand{\cross}{\omega}
\newcommand{\curve}{\Gamma}
\newcommand{\rot}{\mathcal{R}}
\newcommand{\curv}{\mathcal{K}}
\newtheorem{Theorem}{Theorem} 
\newtheorem{Lemma}{Lemma} 
 
\theoremstyle{definition}
\newtheorem{rem}{Remark} 
\newenvironment{Remark}{\begin{rem}}{\qed\end{rem}}
\newtheorem{rems}{Remarks} 

\newtheorem{ass}{Assumption}
\newenvironment{Assumption}{\begin{ass}}{\end{ass}}
\def\OMIT#1{} 
\begin{document}
%
%
\title{{\Large\textbf{
Geometrically induced discrete spectrum in curved tubes 
\footnote{
Accepted for publication in 
\emph{Differential Geometry and its Applications}.
}
}}}
\author{
B.~Chenaud$^1$, 
P.~Duclos$^{2,3}$, 
P.~Freitas$^4$ 
and D.~Krej\v{c}i\v{r}\'{\i}k$^{4}$%
\footnote{
On leave of absence from
\emph{Department of Theoretical Physics,
Nuclear Physics Institute, Academy of Sciences,
250\,68 \v{R}e\v{z} near Prague, Czech Republic}.
} 
}
\date{\footnotesize
\begin{quote}
\emph{
\begin{itemize}
\item[$^1$]
D\'epartement des Sciences Math\'ematiques, Universit\'e Montpellier II, \\ 
Case 051, Place Eug\`ene Bataillon,
34\,095 Montpellier cedex 5, France
\item[$^2$]
Centre de Physique Th\'eorique, CNRS, \\
Luminy, Case 907, 13\,288 Marseille Cedex 9, France 
\item[$^3$]
PhyMat, Universit\'e de Toulon et du Var, \\
BP 132, 83\,957 La Garde, France
\item[$^4$]
Departamento de Matem\'atica, Instituto Superior T\'ecnico, \\
Av. Rovisco Pais, 1049-001 Lisboa, Portugal 
\item[\emph{E-mail:}] 
chenaud@ges.univ-montp2.fr,
duclos@univ-tln.fr, \\
pfreitas@math.ist.utl.pt,
dkrej@math.ist.utl.pt  
\end{itemize}
}
\end{quote}
12 December 2003
}
\maketitle
%
%
\begin{abstract}
\noindent
The Dirichlet Laplacian in curved tubes of arbitrary cross-section 
rotating w.r.t.\ the Tang frame 
along infinite curves in Euclidean spaces of arbitrary dimension 
is investigated. 
If the reference curve is not straight 
and its curvatures vanish at infinity,
we prove that the essential spectrum as a set  
coincides with the spectrum 
of the straight tube of the same cross-section
and that the discrete spectrum is not empty. 
\medskip
\begin{itemize}
\item[\textbf{MSC2000:}] 
81Q10; 58J50; 53A04.
\item[\textbf{Keywords:}] 
Quantum waveguides; Bound states; Dirichlet Laplacian; Tang frame.  
\end{itemize}
\end{abstract}
\newpage
%
%
\section{Introduction}
The relationships between the geometric properties 
of an Euclidean region and the spectrum of the associated 
Dirichlet Laplacian 
constitute one of the classical problems of \emph{spectral geometry},
with important motivations coming 
both from classical and quantum physics.
In this paper we consider this type of interactions
in the case where the region is an infinite tube.
In particular, we are interested in the influence 
of the curvature on the stability of the essential spectrum
and the existence of discrete eigenvalues. 
  
Let~$s \mapsto \curve(s)$ be a unit-speed infinite curve in~$\Real^d$, $d \geq 2$.
Assuming that the curve is $\Smooth^d$-smooth 
and possesses an appropriate $\Smooth^1$-smooth Frenet frame 
$
  \{e_1,\dots,e_d\} 
$
(Assumption~\ref{Ass.Frenet}),
the $i^\mathrm{th}$ curvature~$\kappa_i$ of~$\curve$, $i\in\{1,\dots,d-1\}$,  
is a continuous function of the arc-length parameter $s\in\Real$. 
Given a bounded open connected set~$\cross$ in~$\Real^{d-1}$,  
we define the tube~$\Omega$ of cross-section~$\cross$ about~$\curve$ by
\begin{equation}\label{tube}
  \Omega := \mathcal{L}(\Real\times\cross) 
  \,,
  \qquad
  \mathcal{L}(s,u_2,\dots,u_d)
  :=  \curve(s)+ u_\mu \, \rot_{\mu\nu}(s) \, e_\nu(s)
  \,,
\end{equation}
where $\mu,\nu$ are summation indices taking values in $\{2,\ldots,d\}$
and $(\rot_{\mu\nu})$ is a family of rotation matrices in~$\Real^{d-1}$
chosen in such a way that $(s,u_2,\dots,u_d)$ 
are orthogonal ``coordinates" in~$\Omega$ (\cf~Section~\ref{Sec.rot}),
\ie\/ $\cross$ rotates along~$\curve$ w.r.t.\/ the Tang frame
\cite{Tsao-Gambling_1989}. 
We make the hypotheses (Assumption~\ref{Ass.basic})
that~$\kappa_1$ is bounded,  
$
  \|\kappa_1\|_\infty \sup_{u\in\cross}|u| < 1
$,
and~$\Omega$ does not overlap.

Our object of interest is the Dirichlet Laplacian associated with~$\Omega$,
\ie
\begin{equation}\label{Laplacian}
  -\Delta_D^\Omega
  \qquad\textrm{on}\qquad
  \sii(\Omega).   
\end{equation}
A physical motivation to study this operator for $d=2,3$
comes from the fact that
it is (up to a physical constant) the quantum Hamiltonian
of a free particle constrained to~$\Omega$,
which is widely used to model the dynamics  in mesoscopic systems 
called \emph{quantum waveguides}~\cite{DE,LCM}.  

If~$\curve$ is a straight line (\ie~all $\kappa_i=0$),
then it is easy to see that the spectrum of~(\ref{Laplacian})
is purely absolutely continuous 
and equal to the interval $[\mu_1,\infty)$, 
where
$$
  \mbox{$\mu_1$ is the first eigenvalue 
  of the Dirichlet Laplacian in~$\cross$}
  \,.
$$
The purpose of the present paper is to prove 
that the essential spectrum of the Laplacian~(\ref{Laplacian})
is stable as a set under any curvature which vanishes at infinity, 
and that there is always a geometrically induced spectrum,
\ie\/ the spectrum below~$\mu_1$,
whenever the tube is non-trivially curved.
\begin{Theorem}\label{thm}
Let~$\Omega$ be the infinite tube defined above. 
\begin{itemize}
\item[\emph{(i)}]
If ${\displaystyle \lim_{|s|\to\infty}\kappa_1(s)} = 0$, 
then \
$
  \sigma_\mathrm{ess}(-\Delta_D^\Omega)
  = [\mu_1,\infty) 
  \,;
$
\item[\emph{(ii)}]
If $\kappa_1 \not= 0$, then \
$
  \inf\sigma(-\Delta_D^\Omega)
  < \mu_1 
$.
\end{itemize}
Consequently, if the tube is not straight but it is straight asymptotically,
then $-\Delta_D^\Omega$ has at least one eigenvalue of finite multiplicity
below its essential spectrum, 
\ie~$\sigma_\mathrm{disc}(-\Delta_D^\Omega)\not=\varnothing$.
\end{Theorem}

Here the particularly interesting property is the existence of discrete spectrum,
which is a non-trivial property for unbounded regions. 
From the physical point of view, 
one then deals with quantum \emph{bound states}
of the Hamiltonian~(\ref{Laplacian}), which are known to disturb  
the transport of the particle in the waveguide.

Spectral results of Theorem~\ref{thm} were proved first
by P.~Exner and P.~\v{S}eba~\cite{ES} in 1989 for planar strips (\ie~$d=2$)
under the additional assumptions that the strip width was sufficiently small 
and the curvature~$\kappa_1$ was rapidly decaying at infinity
(roughly speaking as $|s|^{-(\frac{3}{2}+\epsilon)}$).
A few years later, J.~Goldstone and R.~L.~Jaffe~\cite{GJ} proved the results
without the restriction on the strip width,
provided the curvature~$\kappa_1$ had a compact support,
and generalised it to the tubes of circular cross-sections in~$\Real^3$.  
References to other improvements can be found in the review article~\cite{DE},
where the second part of Theorem~\ref{thm} was proved 
under stronger conditions for $d=2,3$ and circular cross-section.
The first part was proved there just for compactly supported~$\kappa_1$
(otherwise, under the additional hypothesis that also the first two derivatives 
of~$\kappa_1$ vanished at infinity,
the authors localised the threshold of the essential spectrum only). 
Let us also mention the paper~\cite{RB}
where a significant weakening of various regularity assumptions 
was achieved for~$d=2$.  
The first part of Theorem~\ref{thm} was proved for $d=2$
in the recent paper~\cite{KKriz}.
The present paper is devoted to a generalisation  
of the results to higher dimensions and cross-sections
rotating along the curve w.r.t. the Tang frame.
 
Our strategy to prove Theorem~\ref{thm} is explained briefly as follows.
Introducing a diffeomorphism from~$\Omega$ 
to the straight tube~$\Omega_0:=\Real \times \cross$ 
by means of the mapping $\mathcal{L}:\Omega_0\to\Omega$,
we transfer the (simple) Laplacian~(\ref{Laplacian}) on (complicated)~$\Omega$
into a unitarily equivalent (complicated) operator~$H$ of the Laplace-Beltrami form 
on (simple)~$\Omega_0$, \cf~(\ref{Laplace.metric}). 
This is the contents of the preliminary Section~\ref{Sec.Preliminaries}.  
The rest of the paper, Section~\ref{Sec.proofs}, 
is devoted to the proof of Theorem~\ref{thm}.
In Section~\ref{Sec.ess}, 
we employ a general characterisation of essential spectrum 
(Lemma~\ref{Iftimie}) adopted from~\cite{DDI} 
in order to establish the first part of Theorem~\ref{thm}. 
The reader will notice that the characterisation is better 
than the classical Weyl criterion in the sense that
it deals with quadratic forms instead of the associated operators themselves, 
\ie\/ we do not need to impose any condition 
on the derivatives of the coefficients of~$H$ in our case.  
The proof of the second part of Theorem~\ref{thm} in Section~\ref{Sec.disc}
is based on the construction of an appropriate trial function for~$H$
inspired by the initial idea of~\cite{GJ}.

Throughout this paper, we use the repeated indices convention
with the range of Latin and Greek indices being
$1,\dots,d$ and $2,\dots,d$, respectively.
The partial derivative w.r.t. a coordinate $x_i$,
$x\equiv(s,u_2,\dots,u_d)\in\Real^d$,
is denoted by a comma with the index~$i$.

\section{Preliminaries}\label{Sec.Preliminaries}
%
\subsection{The reference curve}
Given an integer $d \geq 2$, 
let $\curve:\Real \to \Real^d$ be a unit-speed
$\Smooth^d$-smooth curve satisfying the following hypothesis.
\begin{Assumption}\label{Ass.Frenet}
$\curve$ possesses a positively oriented Frenet frame $\{e_1,\dots,e_d\}$
with the properties that 
\begin{itemize}
\item[(i)]
$
  e_1=\dot{\curve}
$\,;
\item[(ii)]
$
  \forall i \in \{1,\dots,d\},\quad
  e_{i} \in \Smooth^1(\Real,\Real^d)
$\,;
\item[(iii)]
$
  \forall i \in \{1,\dots,d-1\}, \ 
  \forall s\in\Real,
  \quad
  \dot{e}_i(s)
$
lies in the span of $e_1(s),\dots,e_{i+1}(s)$\,.
\end{itemize}
\end{Assumption}
\begin{Remark}
We refer to~\cite[Sec.~1.2]{Kli} for the notion of moving and Frenet frames.
A sufficient condition to ensure the existence
of the Frenet frame of Assumption~\ref{Ass.Frenet}
is to require that for all $s \in \Real$, the vectors
$
  \dot{\curve}(s), \curve^{(2)}(s), \dots, \curve^{(d-1)}(s)
$
are linearly independent, \cf~\cite[Prop.~1.2.2]{Kli}. 
This is always satisfied if $d=2$.
However, we do not assume \emph{a priori}
this non-degeneracy condition for $d \geq 3$
because it excludes the curves such that
$\curve\!\upharpoonright\! I$ lies in a
lower-dimensional subspace of~$\Real^d$
for some open $I\subseteq\Real$. 
We also refer to Remark~\ref{Rem.Exner} below
for further discussions on Assumption~\ref{Ass.Frenet}.
\end{Remark} 

The properties of $\{e_1, \dots, e_d\}$
summarised in Assumption~\ref{Ass.Frenet}
yield the Serret-Frenet formulae, \cf~\cite[Sec.~1.3]{Kli},
\begin{equation}\label{Frenet}
  \dot{e}_i = \curv_{ij} \, e_j
  \,,
\end{equation}
where~$\curv_{ij}$ are coefficients of the skew-symmetric $d \times d$ matrix 
defined by
\begin{equation}\label{curvature}
  (\curv_{ij}) :=
  \begin{pmatrix}
   0                & \kappa_1 &               & \textrm{\Large 0}\\
   -\kappa_1        & \ddots   & \ddots        &                 \\
                    & \ddots   & \ddots        & \kappa_{d-1}    \\
   \textrm{\Large 0} &          & -\kappa_{d-1} & 0
   \end{pmatrix}.
\end{equation}
Here~$\kappa_i:\Real\to\Real$ is called the $i^\mathrm{th}$ curvature of~$\curve$.
Under Assumption~\ref{Ass.Frenet},
the curvatures are continuous functions of the arc-length parameter $s\in\Real$.

\subsection{The Tang frame}\label{Sec.rot}
We introduce now another moving frame along~$\curve$
which better reflects the geometry of the curve
and will be more convenient for our further purposes.

Let the $(d-1)\times(d-1)$ matrix $(\rot_{\mu\nu})$ be defined by
the system of differential equations
\begin{equation}\label{diff.eq}
  \dot{\rot}_{\mu\nu} + \rot_{\mu\rho} \, \curv_{\rho\nu} = 0
\end{equation}
with the initial conditions that $(\rot_{\mu\nu}(s_0))$
is a rotation matrix in $\Real^{d-1}$ for some $s_0 \in \Real$, \ie, 
\begin{equation}\label{rotation}
  \det(\rot_{\mu\nu}(s_0))=1
  \qquad\textrm{and}\qquad
  \rot_{\mu\rho}(s_0) \, \rot_{\nu\rho}(s_0) = \delta_{\mu\nu}.
\end{equation}
Under our assumptions,
the solution of~(\ref{diff.eq}) exists and is continuous
by standard arguments in the theory of differential equations,
\cf~\cite[Sec.~4]{Kurzweil}.
Furthermore, the conditions~(\ref{rotation})
are satisfied for~\emph{all}~$s_0\in\Real$. 
Indeed, by means of Liouville's formula~\cite[Thm.~4.7.1]{Kurzweil} 
and $\tr(\curv_{\mu\nu})=0$, one checks
that $\det\left(\rot_{\mu\nu}\right)=1$ identically, 
while the validity of the second condition for all~$s_0\in\Real$ 
is obtained via the skew-symmetry of~$(\curv_{ij})$:
$$
  (\rot_{\mu\rho} \, \rot_{\nu\rho})^\cdot
  = -\rot_{\mu\rho} \, \rot_{\nu\sigma} \,
  (\curv_{\rho\sigma} + \curv_{\sigma\rho})
  = 0.
$$ 

We set
\begin{equation}
  (\rot_{ij})
  :=
  \begin{pmatrix}
    1 & 0 \\
    0 & (\rot_{\mu\nu})
  \end{pmatrix}
\end{equation}
and define the moving frame
$
  \{\tilde{e}_1, \dots, \tilde{e}_d\}
$
along~$\curve$ by
\begin{equation}\label{Tang}
  \tilde{e}_i := \rot_{ij} \, e_j.
\end{equation}
Combining~(\ref{Frenet}) with~(\ref{diff.eq}) and~(\ref{curvature}), 
one easily finds
\begin{equation}\label{Frenet.bis}
  \dot{\tilde{e}}_1 = \kappa_1 \, e_2
  \qquad\textrm{and}\qquad
  \dot{\tilde{e}}_\mu = - \kappa_1 \, \rot_{\mu 2} \, e_1.
\end{equation}

We call the moving frame 
$
  \{\tilde{e}_1, \dots, \tilde{e}_d\}
$
the \emph{Tang frame} throughout this paper
because it is a natural generalisation
of the Tang frame known from the theory 
of three-dimensional waveguides~\cite{Tsao-Gambling_1989}. 
Its advantage will be clear from the subsequent section.

\subsection{Tubes}
Given a bounded open connected set~$\cross$ in~$\Real^{d-1}$,
let~$\Omega_0$ denote the straight tube~$\Real \times \cross$.
We define the curved tube~$\Omega$ of the same cross-section~$\cross$ 
about~$\curve$ as the image of the mapping, \cf~(\ref{tube}),
\begin{equation}\label{tube.map}
  \mathcal{L}: \Omega_0 \to \Real^d : \
  \left\{
  (s,u) \mapsto \curve(s) + \tilde{e}_\mu(s) \, u_\mu
  \right\},
\end{equation}
\ie, $\Omega:=\mathcal{L}(\Omega_0)$,
where $u \equiv (u_2,\dots,u_d)$. 

Our strategy to deal with the curved geometry of the tube is
to identify~$\Omega$ with the Riemannian manifold~$(\Omega_0,g_{ij})$,
where~$(g_{ij})$ is the metric tensor induced by~$\mathcal{L}$,
\ie~$g_{ij}:=\mathcal{L}_{,i}\cdot\mathcal{L}_{,j}$,
where ``$\cdot$'' denotes the inner product in~$\Real^d$.
(In other words, we parameterise~$\Omega$ globally
by means of the ``coordinates''~$(s,u)$.)
To this aim, we need to impose a natural restriction on~$\Omega$ 
in order to ensure that
$\mathcal{L}:\Omega_0\to\Omega$ is a $\Smooth^1$-diffeomorphism.
Namely, defining
\begin{equation}\label{a}
  a:=\sup_{u\in\cross}|u| \,,
\end{equation}
where $|u|:=\sqrt{u_\mu u_\mu}$, we make the hypothesis
\begin{Assumption}\label{Ass.basic} 
\rule{1ex}{0ex}
\begin{itemize}
\item[(i)]
$\kappa_1\in\sinf(\Real)$ \ and \ $a\,\|\kappa_1\|_\infty < 1$\,; 
\item[(ii)]
$\Omega$ does not overlap\,.
\end{itemize}
\end{Assumption}

Using formulae~(\ref{Frenet.bis}), one easily finds
\begin{equation}\label{metric}
  (g_{ij}) = \diag\left(h^2,1,\dots,1\right)
  \qquad\textrm{with}\qquad
  h(s,u) := 1-\kappa_1(s) \, \rot_{\mu 2}(s) \, u_\mu \,.
\end{equation}
By virtue of the inverse function theorem, 
the mapping $\mathcal{L}:\Omega_0\to\Omega$ is a local $\Smooth^1$-diffeomorphism 
provided~$h$ does not vanish on~$\Omega_0$,
which is guaranteed by the condition~(i) of Assumption~\ref{Ass.basic}
because
\begin{equation}\label{C+-}
  0 < C_- \leq h(s,u) \leq C_+ < 2
  \qquad\textrm{with}\quad
  C_\pm := 1 \pm a\,\|\kappa_1\|_\infty,
\end{equation}
where we have used that 
$\rot_{\mu 2}\,\rot_{\mu 2}=1$ by~(\ref{rotation})
and $\sqrt{u_\mu u_\mu}<a$ by~(\ref{a}). 
The mapping then becomes a global diffeomorphism if 
it is required to be injective in addition, 
\cf~the condition~(ii) of Assumption~\ref{Ass.basic}.
\begin{Remark}\label{Rem.basic} 
Formally, it is possible to consider $(\Omega_0,g_{ij})$
as an abstract Riemannian manifold where only the curve~$\curve$
is embedded in~$\Real^d$.
Then we do not need to assume the condition~(ii) of Assumption~\ref{Ass.basic}.
\end{Remark}

Note that the metric tensor~(\ref{metric}) is diagonal
due to our special choice of the ``transverse'' frame
$\{\tilde{e}_2,\dots,\tilde{e}_d\}$,
which is the advantage of the Tang frame.
At the same time, it should be stressed here that
while the shape of the tube~$\Omega$ is not influenced
by a special choice of the rotation~$(\rot_{\mu\nu})$
provided~$\cross$ is circular,
this may not be longer true for a general cross-section.  
In this paper, we choose rotations determined 
by the Tang frame due to the technical simplicity.

We set $g := \det(g_{ij}) = h^2$, which defines through
$d\vol := h(s,u)\,ds\,du$ the volume element of~$\Omega$;
here~$du$ denotes the $(d-1)$-dimensional Lebesgue measure in~$\cross$.
\begin{Remark}[Low-dimensional examples]\label{Rem.Examples.tilde}
When $d=2$,
the cross-section~$\cross$ is just the interval~$(-a,a)$,
the curve~$\curve$ has only one curvature~$\kappa_1 =: \kappa$,
the rotation matrix $(\rot_{\mu\nu})$ equals (the number) $1$
and
$$
  h(s,u)=1-\kappa(s)\,u.
$$
When $d=3$, it is convenient to make the Ansatz
$$
  (\rot_{\mu\nu}) =
  \begin{pmatrix}
    \cos\theta & -\sin\theta \\
    \sin\theta & \cos\theta
  \end{pmatrix},
$$
where~$\theta$ is a real-valued differentiable function.
Then it is easy to see that~(\ref{diff.eq}) reduces
to the differential equation $\dot{\theta}=\tau$,
where~$\tau$ is the torsion of~$\curve$,
\ie~one put $\kappa:=\kappa_1$ and $\tau:=\kappa_2$.
Choosing~$\theta$ as an integral of~$\tau$, we can write
$$
  h(s,u) =
  1-\kappa(s)\left(\cos\theta(s)\,u_2+\sin\theta(s)\,u_3\right).
$$
\end{Remark}
\begin{Remark}[On Assumption~\ref{Ass.Frenet}]\label{Rem.Exner}
As pointed out by P.~Exner~\cite{E-private},
the existence of the Frenet frame required in Assumption~\ref{Ass.Frenet}  
is rather a technical hypothesis only.
Indeed, what we actually need is that $\mathcal{L}:\Omega_0\to\Omega$,
with~$\mathcal{L}$ given by~(\ref{tube.map}), is a $\Smooth^1$-diffeomorphism, 
and this is possible to ensure in certain situations
even if Assumption~\ref{Ass.Frenet} does not hold.  
To see this, let~$\cross$ be circular 
and consider a curve~$\curve$ possessing
the required Frenet frame on $(-\infty,0)$ and $(0,\infty)$,  
$e_1\in\Smooth^1(\{0\})$, but $e_\mu\not\in\Smooth^0(\{0\})$
in the sense that 
$
  e_\mu(0+) = S_{\mu\nu} e_\nu(0-)
$
for each $\mu\in\{2,\dots,d\}$,
where $(S_{\mu\nu})\not=1$ is a constant matrix satisfying relations
analogous to~(\ref{rotation}),
\ie\/ the transverse frames $\{e_2(0+),\dots,e_d(0+)\}$ 
and $\{e_2(0-),\dots,e_d(0-)\}$ are rotated to each other
(see~\cite[Chap.~1, pp.~34]{Spivak2} for an example of such a curve in~$\Real^3$). 
Since the rotation matrix~$(R_{\mu\nu})$ is determined
uniquely up to a multiplication by a constant rotation matrix,
the Tang frame defined by~(\ref{Tang}) can be chosen to be continuous at zero
by the requirement $R_{\mu\nu}(0+)S_{\nu\rho}=R_{\mu\rho}(0-)$.
The $\Smooth^1$-continuity at zero then follows by~(\ref{Frenet.bis}) 
together with the fact that necessarily $\kappa_1(0)=0$.
\end{Remark}
%

\subsection{The Laplacian}\label{Sec.Laplacian}
Our strategy to investigate the Laplacian~(\ref{Laplacian})
is to express it in the coordinates determined by~(\ref{tube.map}). 
More specifically, using the mapping~(\ref{tube.map}),
we identify the Hilbert space $\sii(\Omega)$
with $\sii(\Omega_0,d\vol)$
and consider on the latter the sesquilinear form
\begin{equation}\label{form}
  Q(\psi,\phi) :=
  \int_{\Omega_0} \overline{\psi_{,i}}\,g^{ij}\,\phi_{,j}
  \ d\vol \,,
  \qquad
  \psi,\phi\in\Dom Q :=
  \sobi(\Omega_0,d\vol),
\end{equation}
where~$g^{ij}$ denotes the coefficients 
of the inverse of the metric tensor~(\ref{metric}).
The form~$Q$ is clearly densely defined,
non-negative, symmetric and closed on its domain.
Consequently, there exists a non-negative self-adjoint
operator~$H$ associated with~$Q$
which satisfies~$\Dom H \subset \Dom Q$.
We have
\begin{align}\label{Laplace.metric} 
  \Dom H &=
  \left\{\psi\in\sobi(\Omega_0,d\vol) \left| \
  \partial_i g^\demi g^{ij} \partial_j \psi
  \in\sii(\Omega_0,d\vol) 
  \right.\right\},
  \nonumber \\
  \forall\psi\in\Dom H, \quad
  H\psi&=-g^{-\demi}\partial_i \big(g^\demi g^{ij} \partial_j \psi\big).
\end{align}
Actually, (\ref{Laplace.metric}) is a general expression 
for the Laplace-Beltrami operator 
in a manifold equipped with a metric~$(g_{ij})$.
Using the particular form~(\ref{metric}) of our metric tensor, 
we can write
\begin{equation}\label{Laplace.metric.h}
  H = - \frac{1}{h} \, \partial_1 \, \frac{1}{h} \, \partial_1
  - \partial_\mu \partial_\mu
  + \frac{\kappa_1\,\rot_{\mu 2}}{h} \, \partial_\mu
\end{equation}
in the form sense.

The norm and the inner product in the Hilbert space $\sii(\Omega_0,d\vol)$
will be denoted by $\|\cdot\|_g$ and $(\cdot,\cdot)_g$, respectively.
The usual notation without the subscript~``$g$''
will be reserved for the similar objects in $\sii(\Omega_0)$.

\begin{Remark}[Unitarily equivalent operator]
Assuming that the reference curve~$\curve$ is $\Smooth^{d+1}$-smooth, 
it is possible to ``rewrite''~$H$ into a Schr\"odinger-type operator
acting on the Hilbert space~$\sii(\Omega_0)$,
without the additional weight~$g^\demi$ in the measure of integration. 
Indeed, defining $\hat{H}:=U H U^{-1}$, where
$
  U: \psi \mapsto g^\frac{1}{4}\psi
$
is a unitary transformation
from $\sii(\Omega_0,d\vol)$ to~$\sii(\Omega_0)$,
we get 
$$
  \hat{H} =  
  -g^{-\frac{1}{4}}\partial_i g^\demi g^{ij}
  \partial_j g^{-\frac{1}{4}}
  \qquad\textrm{on}\quad\
  \sii(\Omega_0)
$$
in the form sense. 
Commuting then $g^{-\frac{1}{4}}$ with the gradient components,
we can write
\begin{equation}\label{Hamiltonian.potential}
  \hat{H} = -\partial_1 \, \frac{1}{h^2} \, \partial_1
  + \partial_\mu \partial_\mu + V
\end{equation}
in the form sense, where
\begin{align} 
  V \ := \ &
  - \frac{5}{4} \, \frac{(h_{,1})^2}{h^4}
  + \demi \, \frac{h_{,11}}{h^3}
  - \frac{1}{4} \, \frac{h_{,\mu}\,h_{,\mu}}{h^2}
  + \demi \, \frac{h_{,\mu\mu}}{h}
  \label{potential} \\
  \ = \ &
  - \frac{1}{4} \, \frac{\kappa_1^2}{h^2}
  + \demi \, \frac{h_{,11}}{h^3}
  - \frac{5}{4} \, \frac{(h_{,1})^2}{h^4}
  \,.
  \label{potential.bis}
\end{align}
Actually, (\ref{Hamiltonian.potential}) with~(\ref{potential})
is a general formula valid for any $\Smooth^1$-smooth metric
of the form $(g_{ij})=\diag(h^2,1,\dots,1)$.
(Note that the required regularity is indeed sufficient 
if the formula for the potential~(\ref{potential}) 
is understood in the weak sense of forms.)
In our special case when~$h$ is given by~(\ref{metric}),
we find easily that
$h_{,\mu}(\cdot,u)=-\kappa_1 \, \rot_{\mu 2}$,
$h_{,\mu\nu}=0$,
and~(\ref{potential.bis}) follows at once.
Moreover, (\ref{diff.eq}) gives
\begin{align*} 
  h_{,1}(\cdot,u)&=u_\mu\,\rot_{\mu\alpha}
  \big(\dot{\curv}_{\alpha 1}
  -\curv_{\alpha\beta}\curv_{\beta 1}\big)\,,
  \\
  h_{,11}(\cdot,u)&=u_\mu\,\rot_{\mu\alpha}
  \big(\ddot{\curv}_{\alpha 1}
  -\dot{\curv}_{\alpha\beta}\curv_{\beta 1}
  -2\;\!\curv_{\alpha\beta}\dot{\curv}_{\beta 1}
  +\curv_{\alpha\beta}\curv_{\beta\gamma}\curv_{\gamma 1}
  \big)
  \,. 
\end{align*}
\end{Remark}
We shall neither need nor use the unitarily equivalent operator
from the above remark, however, for motivation purposes,  
it is interesting to notice that the potential~$V$ becomes attractive
if~$a$ is sufficiently small and the curvatures, 
together with some of their derivatives, vanish at infinity.
Since the latter also implies that~$h$ tends to~$1$ as $a \to 0$,
it is easy to see that~$\hat{H}$ has always discrete eigenvalues 
below its essential spectrum for~$a$ small enough.
In this paper, we prove this property under an asymptotic condition
which involves the curvature~$\kappa_1$ only (\cf~(\ref{decay}) below)
and without any restriction on~$a$
(except for the natural one in Assumption~\ref{Ass.basic}).     
We also note that various techniques from the theory of Schr\"odinger operators
can be applied to~$\hat{H}$, \cf~\cite{DE}.   

\subsection{Straight tubes}\label{Sec.Straight}
If the tube is straight in the sense that each $\kappa_i = 0$,  
then the Laplacian~(\ref{Laplacian}) coincides with the decoupled operator
\begin{equation}\label{straight.tube}
  H_0
  := \overline{
  -\Delta^\Real \otimes \Id
  + \Id\otimes(-\Delta_D^\cross) }
  \qquad\textrm{on}\quad\
  \sii(\Real)\otimes\sii(\cross),
\end{equation}
where~$\Id$ denotes the identity operator on the corresponding spaces
and the bar stands for the closure.
The operators $-\Delta^\Real$ and $-\Delta_D^\cross$
denote the usual Laplacian on $\sii(\Real)$
and the Dirichlet Laplacian on $\sii(\cross)$, respectively.
Alternatively, $H_0$ can be introduced as the operator associated 
with the form~$Q_0$
defined by~(\ref{form}), where now the metric tensor 
is the identity matrix~$(\delta_{ij})$ 
and $d\vol = ds\,du$ is the Lebesgue measure in~$\Real\times\cross$.

The operator $-\Delta_D^\cross$ has a purely discrete spectrum 
consisting of eigenvalues 
$
  \mu_1 < \mu_2 \leq \dots \mu_n < \dots;
$ 
the corresponding eigenfunctions are denoted as~$\ef_n$
and we normalise them in such a way that 
$
  \|\ef_n\|_{\sii(\cross)}=1 .
$
The lowest eigenvalue~$\mu_1$ is, of course, positive, simple
and the eigenfunction~$\ef_1$ can be chosen positive.  

In view of the decomposition~(\ref{straight.tube}), 
\begin{equation}\label{StraightSpectrum}
  \sigma(H_0)
  =\sigma_\mathrm{ess}(H_0)
  =[\mu_1,\infty) 
\end{equation}
and the spectrum is absolutely continuous.

\section{Proofs}\label{Sec.proofs}
%
\subsection{The essential spectrum}\label{Sec.ess}
We prove that the essential spectrum of a curved tube~$\Omega$ 
coincides with the one of~$\Omega_0$ provided 
the former is straight asymptotically in the sense that 
\begin{equation}\label{decay}
  \lim_{|s|\to\infty} \kappa_1(s) = 0 .
\end{equation}
Our method is based on the following characterisation 
of the essential spectrum of~$H$.
\begin{Lemma}\label{Iftimie}
$\lambda\in\sigma_\mathrm{ess}(H)$ if and only if
there exists $\{\psi_n\}_{n\in\Nat}\subset\Dom Q$ 
such that
\begin{itemize}
\item[\emph{(i)}] \quad 
$
  \forall n\in\Nat,\quad \|\psi_n\|_g = 1
$,
\item[\emph{(ii)}] \quad 
$
  \forall n\in\Nat,\quad 
  \supp\psi_n \subset \left\{(s,u)\in\Omega_0 \,| \ |s| \geq n \right\}
$, 
\item[\emph{(iii)}] \quad
  $(H-\lambda)\psi_n \xrightarrow[n\to\infty]{ } 0
  \quad\textrm{in}\ \left(\Dom Q\right)^*$.
\end{itemize}
\end{Lemma}
\noindent
Here $(\Dom Q)^*$ denotes the dual of the space $\Dom Q$ defined in~(\ref{form}).
We note that $H+1:\Dom Q\to(\Dom Q)^*$ is an isomorphism and 
\begin{equation}\label{-1norm}
  \|\psi\|_{-1,g} 
  := \|\psi\|_{(\Dom Q)^*} 
  = \sup_{\phi\in (\Dom Q)\setminus\{0\}}
  \frac{|(\phi,\psi)_g|}{\|\phi\|_{1,g}} 
\end{equation}
with
$$
  \|\phi\|_{1,g} 
  := \|\phi\|_{\Dom Q}
  = \sqrt{Q[\phi]+\|\phi\|_g^2} \,.
$$

The proof of the above lemma is quite similar 
to the proof of Lemma~4.2 in~\cite{DDI}.
It is based on a general characterisation of essential spectrum,
\cite[Lemma~4.1]{DDI}, which is better than the Weyl criterion  
in the sense that the former requires to find 
a sequence from the form domain of~$H$ only
(\cf~the statement of Lemma~\ref{Iftimie} and the required property~(iii)
with the Weyl criterion~\cite[Thm.~7.24]{Weidmann}).
The second property~(ii) reflects the fact that the essential spectrum 
is determined by the geometry at infinity only.
\begin{Remark}\label{Rem.equivalence}
Since the metric~$(g_{ij})$ is uniformly elliptic due to~(\ref{C+-}), 
the norms in the spaces $\sii(\Omega_0,d\vol)$ and $\sobi(\Omega_0,d\vol)$ 
are equivalent with those of $\sii(\Omega_0)$ and $\sobi(\Omega_0)$, 
respectively, and the respective spaces can be identified as sets. 
In particular, 
$$
  C_- \|\psi\|^2 \leq \|\psi\|_g^2 \leq C_+ \|\psi\|^2 ,
$$
and similarly for $\|\cdot\|_{1,g}$ and $\|\cdot\|_{-1,g}$.
\end{Remark}
\paragraph{Proof of Theorem~\ref{thm}, part~(i).}
Let $\lambda\in\sigma_\mathrm{ess}(H_0) \equiv [\mu_1,\infty)$. 
By Lemma~\ref{Iftimie}, there exists a sequence 
$\{\tilde{\psi}_n\}_{n\in\Nat} \subset \Dom Q_0$ 
such that it satisfies the properties (i)--(iii) of the Lemma
for $g_{ij}=\delta_{ij}$, $H=H_0$ and $Q=Q_0$. 
We will show that the sequence $\{\psi_n\}_{n\in\Nat}$
defined by $\psi_n:=\tilde{\psi}_n/\|\tilde{\psi}_n\|_g$ 
for every $n\in\Nat$
satisfies the properties (i)--(iii) of Lemma~\ref{Iftimie}
for $g_{ij}$, $H$ and $Q$, 
\ie~$\sigma_\mathrm{ess}(H_0) \subseteq \sigma_\mathrm{ess}(H)$.
First of all, notice that $\psi_n$ is well defined and belongs to $\Dom Q$
for every $n\in\Nat$ due to Remark~\ref{Rem.equivalence}.
Moreover, writing
$$
  \tilde{\psi}_n 
  = (1+H_0)^{-1}(H_0-\lambda)\tilde{\psi}_n 
  + (1+H_0)^{-1} (\lambda+1) \tilde{\psi}_n
  \,, 
$$
we see that the sequence $\{\tilde{\psi}_n\}_{n\in\Nat}$ is bounded in $\Dom Q_0$,
and therefore $\{\psi_n\}_{n\in\Nat}$ is bounded in $\Dom Q$
by Remark~\ref{Rem.equivalence}.
Since the conditions~(i) and~(ii) hold trivially true for~$\{\psi_n\}_{n\in\Nat}$,
it remains to check the third one.
By the definition of~$H$, \cf~(\ref{form}),  
we can write
\begin{multline*}
  \big(\phi,(H-\lambda)\psi_n\big)_g
  \\ 
  = \big(\phi,(H_0-\lambda)\psi_n\big)
  + \big(\phi_{,i},(g^\demi g^{ij}-\delta^{ij})\psi_{n,j}\big)
  - \lambda \big(\phi,(g^\demi-1)\psi_n\big)
\end{multline*}
for every $\phi\in\Dom Q$.
The Minkowski inequality, the formula~(\ref{-1norm}), 
the fact that $(H_0-\lambda)\tilde{\psi}_n \to 0$ in $(\Dom Q_0)^*$ as $n\to\infty$, 
and a repeated use of Remark~\ref{Rem.equivalence} 
yield that it is enough to show that  
$$
  \sup_{\phi\in (\Dom Q)\setminus\{0\}}
  \frac{\big| 
  \big(\phi_{,i},(g^{ij}-g^{-\demi} \delta^{ij})\psi_{n,j}\big)_g 
  \big|
  + \lambda \big| \big(\phi,(1-g^{-\demi})\psi_n\big)_g
  \big|}
  {\|\phi\|_{1,g}} 
  \xrightarrow[n\to\infty]{} 0.
$$
However, the latter is easily established by means of the Schwarz inequality, 
the estimates $\|\phi_{,i}\|_g, \|\phi\|_g \leq \|\phi\|_{1,g}$,
the fact that $\{\psi_n\}_{n\in\Nat}$ is bounded in $\Dom Q$,
and the expression for the metric~(\ref{metric})
together with~(\ref{decay}) and the property~(ii) of Lemma~\ref{Iftimie}. 

One proves that $\sigma_\mathrm{ess}(H) \subseteq \sigma_\mathrm{ess}(H_0)$
in the same way.
\qed
\begin{Remark}
It is clear from the previous proof that
a stronger result than the first part of Theorem~\ref{thm} can be proved.
If~$h$ and~$\tilde{h}$ are two positive functions 
(determining through~(\ref{metric})
two tube metrics~$(g_{ij})$ and~$(\tilde{g}_{ij})$, respectively)
such that 
$
  \sup_{u\in\cross}|h(s,u)-\tilde{h}(s,u)| \to 0 
$  
as~$|s|\to\infty$,
then the essential spectra of the corresponding
operators~$H$ and~$\tilde{H}$ (given by~(\ref{Laplace.metric}) 
with~$(g_{ij})$ and~$(\tilde{g}_{ij})$, respectively) coincide as sets.
\end{Remark}

Let us finally notice that a detailed study of the \emph{nature}
of the essential spectrum in curved tubes
has been performed in~\cite{KT};
in particular, the absence of singular continuous spectrum
is proved there under suitable assumptions 
about the decay of curvature at infinity. 

\subsection{The geometrically induced spectrum}\label{Sec.disc}
In this section we show that $\inf\sigma(H) < \mu_1$
whenever $\kappa_1 \not= 0$,
\ie~there is always a spectrum below the energy~$\mu_1$
in non-trivially curved tubes~$\Omega$.
We call it geometrically induced spectrum
because it does not exist for the straight tube~$\Omega_0$, 
\cf~(\ref{StraightSpectrum}).
Furthermore, it follows by the part~(i) of Theorem~\ref{thm}
that this geometrically induced spectrum is discrete
if we suppose~(\ref{decay}) in addition.

Our proof is based on the variational strategy of finding
a trial function~$\Psi$ from the form domain of~$H$ such that
\begin{equation}
  Q_1[\Psi] := Q[\Psi] - \mu_1 \, \|\Psi\|_g^2 < 0 . 
\end{equation}
The construction of such a~$\Psi$ follows 
the initial idea of~\cite{GJ}
and the subsequent improvements of~\cite{RB}
and~\cite[Thm.~2.1]{DE}.  
\paragraph{Proof of Theorem~\ref{thm}, part~(ii).}
Let $\{\Psi_n\}_{n\in\Nat}\subset\Dom Q$ and $\Phi \in \Dom Q$.
Defining  
$
  \Psi_{n,\varepsilon} := \Psi_n + \varepsilon\,\Phi  
$
for every $(n,\varepsilon)\in\Nat\times\Real$,
we can write
$$
  Q_1[\Psi_{n,\varepsilon}]
  = Q_1[\Psi_n] + 2\varepsilon \, Q_1(\Phi,\Psi_n) + \varepsilon^2 Q_1[\Phi].
$$
Our strategy will be to choose $\{\Psi_n\}_{n\in\Nat}$ 
and~$\Phi$ so that 
\begin{equation}\label{DE.trick}
  \lim_{n\to\infty} Q_1[\Psi_n] = 0
  \qquad\textrm{and}\qquad
  \lim_{n\to\infty} Q_1(\Phi,\Psi_n) \not = 0.
\end{equation}
Then we can choose a sufficiently large $n\in\Nat$  
and a sufficiently small $\varepsilon\in\Real$ with a suitable sign
so that $Q_1[\Psi_{n,\varepsilon}] < 0$, which proves the claim.
 
We put $\Psi_n:=\varphi_n\otimes\ef_1$, 
where~$\ef_1$ is the first eigenfunction of~$-\Delta_D^\cross$,
\cf~Section~\ref{Sec.Straight}, 
and $\{\varphi_n\}_{n\in\Nat}$ is a mollifier of~$1$ in $\Sobi(\Real)$,
\ie~a family of functions~$\varphi_n$ from~$\Sobi(\Real)$ satisfying:
\begin{itemize}
\item[(i)]
$
  \forall n\in\Nat, \quad 0 \leq \varphi_n \leq 1 
$,
\item[(ii)]
$
  \varphi_n(s) \xrightarrow[n\to\infty]{} 1
  \quad \textrm{for a.e.} \ s\in\Real
$,
\item[(iii)]
$
  \|\dot{\varphi}_n\|_{\sii(\Real)} \xrightarrow[n\to\infty]{} 0
$.
\end{itemize}
(Probably the simplest example of such a family is given 
by the continuous even $\varphi_n$'s 
such that they are equal to~1 on $[0,n)$,
with a constant derivative on $[n,2n+1)$,
and equal to~$0$ on $[2n+1,\infty)$.)
Using the expression~(\ref{Laplace.metric.h}) for the Laplacian
and the fact that $(\partial_\mu\partial_\mu+\mu_1) \ef_1=0$,  
we obtain immediately that
$$
  Q_1[\Psi_n] 
  = \left(\Psi_{n,1},h^{-1}\Psi_{n,1}\right)
  + \left(\Psi_n,\kappa_1\rot_{\mu2}\Psi_{n,\mu}\right).
$$
The second term at the r.h.s. is equal to zero by an integration by parts,
while the first (positive) one can be estimated from above by
$
  C_-^{-1} \|\Psi_{n,1}\|^2 
$
due to~(\ref{C+-}). Since 
$
  \|\Psi_{n,1}\| = \|\dot{\varphi}_n\|_{\sii(\Real)}
$
by the normalisation of~$\ef_1$,
we verify the first property of~(\ref{DE.trick}).  

The second property is checked if we take, for instance,
$$
  \Phi(s,u) := \phi(s) \, \rot_{\mu2}(s) \, u_\mu \, \ef_1(u) 
  \ \in\Dom Q ,
$$
where~$\phi\in\Sobi(\Real)\setminus\{0\}$ is a non-negative function   
with a compact support contained in an interval 
where~$\kappa_1$ is not zero and does not change sign
(such an interval surely exists because~$\kappa_1\not=0$ is continuous).  
Indeed, in the same way as above, we find  
$$
  Q_1(\Phi,\Psi_n)
  = \left(\Phi_{,1},h^{-1}\Psi_{n,1}\right)
  + \left(\Phi,\kappa_1\rot_{\mu2}\Psi_{n,\mu}\right) ,
$$
where the first term at the r.h.s. tends to zero as~$n\to\infty$ 
because its absolute value can be estimated by
$
  C_-^{-1} \|\Phi_{,1}\| \|\Psi_{n,1}\|
$,
while the second one is equal to
$$
  - \mbox{$\demi$} \big(\phi\,\ef_1,
  \kappa_1 \rot_{\mu 2}\rot_{\mu 2}\,\varphi_n\,\ef_1\big)
  = - \mbox{$\demi$} \left(\phi,\kappa_1\varphi_n\right)_{\sii(\Real)}
$$
by an integration by parts; 
the last identity then holds due to~(\ref{rotation}) 
and the normalisation of~$\ef_1$.
Summing up, we conclude that 
$$
  \lim_{n\to\infty} Q_1(\Phi,\Psi_n) 
  = - \demi \int_{\supp\phi} \phi(s)\,\kappa_1(s)\,ds
  \not= 0
$$
by the dominated convergence theorem. 
\qed
\begin{Remark}
Suppose Assumptions~\ref{Ass.Frenet} and~\ref{Ass.basic}.
If the tube~$\Omega$ is non-trivially curved and asymptotically straight,
it follows by Theorem~\ref{thm} that 
$
  \sigma_\mathrm{disc}(-\Delta_D^\Omega)
  \subset[0,\mu_1)
$
and it is not empty. 
Furthermore, it can be shown by standard arguments
(see, \eg, \cite[Sec.~8.12]{Gilbarg-Trudinger})
that the minimum eigenvalue, \ie~$\inf\sigma(-\Delta_D^\Omega)$,
is simple and has a positive eigenfunction. 
One also has $\inf\sigma(-\Delta_D^\Omega)>0$.
(Actually, a stronger lower bound to the spectral threshold
has been derived in~\cite{EFK}.)  
\end{Remark}
%

\section*{Acknowledgements}
\addcontentsline{toc}{section}{Acknowledgements}
The authors wish to thank Pavel Exner for useful discussions
and Remark~\ref{Rem.Exner}. 
This work was partially supported by FCT/POCTI/FEDER, Portugal,
and GA AS CR grant IAA 1048101.

%
%
\addcontentsline{toc}{section}{References}

\providecommand{\bysame}{\leavevmode\hbox to3em{\hrulefill}\thinspace}
\providecommand{\MR}{\relax\ifhmode\unskip\space\fi MR }
\providecommand{\MRhref}[2]{%
  \href{http://www.ams.org/mathscinet-getitem?mr=#1}{#2}
}
\providecommand{\href}[2]{#2}

\end{document}